\documentclass[reqno, a4paper, 10pt]{amsart}
\usepackage{amsmath,amsthm,amssymb,amscd,amstext,amsopn,amsxtra,amsfonts, bbm, mathrsfs}
\usepackage{verbatim}

\usepackage[all]{xy}
\SelectTips{eu}{}
\entrymodifiers={+!!<0pt,\fontdimen22\textfont2>}

\usepackage[pdftex]{hyperref}
\usepackage{url}

\swapnumbers
\theoremstyle{plain}
\newtheorem{thm}[subsection]{Theorem}
\newtheorem{lemma}[subsection]{Lemma}
\newtheorem{prop}[subsection]{Proposition}

\theoremstyle{remark}

\theoremstyle{definition}

\newtheorem{remark}[subsection]{Remark}

\numberwithin{equation}{subsection}

\def\cH{\mathcal{H}}

\def\cO{\mathcal{O}}

\def\11{\mathbf{1}}

\def\DD{\mathbf{D}}

\def\ZZ{\mathbf{Z}}

\def\dim{\mathrm{dim\,}} 

\def\Ext{\mathrm{Ext}}

\def\Flag{\mathcal{B}}

\title{A remark on some bases in the Hecke algebra}
\author{R. Virk}
\address{Department of Mathematics\\
University of California\\
Davis, CA 95616}
\email{virk@math.ucdavis.edu}

\begin{document}

\maketitle
\setcounter{tocdepth}{1}
\tableofcontents
\section{Introduction}
In this brief note we consider some bases in the Hecke algebra and exhibit certain dualities between them. Our results and arguments are completely combinatorial. However, to explain the motivation and put the contents of this note in perspective we need to discuss some results from geometric representation theory.

Let $G$ be a reductive group, $B\subseteq G$ a Borel subgroup and let $\Flag$ be the flag variety of $G$. Then the Iwahori-Hecke algebra of $G$ realizes the Grothendieck group of $B$-equivariant perverse sheaves on $\Flag$. Multiplication in the Hecke algebra corresponds to convolution of sheaves.

In his streamlined treatment of the Hecke algebra, Soergel \cite{So97} considers several pairwise commuting automorphisms (see \cite[Thm.\ 2.7]{So97}. The geometric significance of some of these is well known. For instance, the automorphism $d$ \eqref{defd} corresponds to Verdier duality on sheaves. This note was motivated by trying to understand the geometric significance of the automorphism $b$ \eqref{defb}. I will now outline the contents of this note from this perspective.

A cornerstone of Kazhdan-Lusztig theory is 
the identification of IC-complexes (=simple objects) on $\Flag$ with the Kazhdan-Lusztig basis of the Hecke algebra. In Thm.\ \ref{projectivesexist} we identify/construct the basis corresponding to projective objects. In Thm.\ \ref{tiltingduality} we explain how $b$ connects the projective basis with the Kazhdan-Lusztig basis. Actually, Thm.\ \ref{tiltingduality} implies that $b$ switches the Kazhdan-Lusztig basis with the basis corresponding to tilting objects (see \cite[\S9]{Virk} for the precise statement). Now, Soergel's `tilting duality' \cite{So} switches tilting objects with projective objects and the Koszul duality of \cite{BGS} switches projective objects with simple objects. 
It follows that $b$ is the composition of these two dualities. Another way to see this is to use \cite[Thm. 4.4]{So3} to directly deduce that $b$ switches the Kazhdan-Lusztig basis with the tilting basis. The aforementioned Koszul duality does not commute with convolution (it doesn't preserve the monoidal unit). In particular, it does not descend to a ring automorphism of the Hecke algebra. However, Thm.\ \ref{tiltingduality} is a combinatorial shadow of the statement that the tilting duality intertwines convolution and Koszul duality (cf.\ \cite[Conjecture 5.18, Thm.\ 5.24, Thm.\ 6.10]{BeGi}). \\

\noindent\textbf{Acknowledgments} I am grateful to W.\ Soergel for some very helpful correspondence and to Arun Ram for his comments on a preliminary version of this note. This work is supported by the NSF grant DMS-0652641.

\section{The Hecke algebra}\label{s:hecke}
\subsection{}
Let $(W,S)$ be a Coxeter system, $\ell\colon W\to \ZZ_{\geq 0}$ the corresponding length function and let $\leq$ denote the Bruhat order on $W$. In particular, $x<y$ means $x\leq y$ and $x\neq y$. The identity element in $W$ is denoted by $e$.

\subsection{}
The Hecke algebra $\cH$ is the free $\ZZ[v, v^{-1}]$-module $\bigoplus_{x\in W}\ZZ[v, v^{-1}]T_x$ with $\ZZ[v, v^{-1}]$-algebra structure given by
\begin{align}
\label{braid}T_xT_y &= T_{xy} &&\mbox{if $\ell(xy)=\ell(x)+\ell(y)$}, \\
\label{quadratic}(T_s+1)(T_s-v^{-2}) &= 0 &&\mbox{for all $s\in S$}.
\end{align}
Set $H_x = v^{\ell(x)}T_x$ for all $x\in W$. Then \eqref{quadratic} implies that each $H_x$ is invertible. In particular, if $s\in S$, then $H_s^2 = 1 - (v-v^{-1})H_s$. Hence,
\begin{equation}\label{inverse}
H_s^{-1} = H_s + (v-v^{-1}).
\end{equation}
Note that $T_e=1$.

\subsection{}
Define a ring involution $d\colon\cH\to\cH$ by
\begin{equation}\label{defd}d(v)=v^{-1}, \qquad d(H_x) =H_{x^{-1}}^{-1}.\end{equation}
We will often write $\overline{H}$ instead of $d(H)$.
We call $C\in\cH$ self-dual if $\overline{C} = C$. For each $s\in S$ set
\begin{equation}
C_s= H_s + v.
\end{equation}
Then, by \eqref{inverse}, each $C_s$ is self-dual. In the sequel we will need the following well known formulae:
\begin{lemma}\label{translationthroughthewall}
Let $s\in S$ and let $x\in W$ be arbitrary. Then
\[ C_sH_x = \begin{cases}
H_{sx} + vH_x & \mbox{if $sx>x$}; \\
H_{sx}+v^{-1}H_x & \mbox{if $sx<x$}.
\end{cases}\]
\end{lemma}
\begin{proof}
This is a direct computation and is left to the reader.
\end{proof}

The following classical result is due to Kazhdan and Lusztig \cite{KL79}, the proof presented here is stolen from \cite{So97}.
\begin{lemma}\label{klbasisstrong}For each $x\in W$ there exists a self dual element $C_x\in\cH$ such that $C_x \in H_x + \sum_{y<x} v\ZZ[v]H_y$.
\end{lemma}

\begin{proof}
Proceed by induction on the Bruhat order. Certainly we can start our induction with $C_e=H_e=1$. Now let $x\in W$ be given and suppose we know the existence of $C_y$ for all $y<x$. If $x\neq e$ find $s\in S$ such that $sx<x$. Then by the induction hypothesis
$C_s C_{sx} = H_x + \sum_{y<x}h_y H_y$
for some $h_y\in \ZZ[v]$. Set
$C_x = C_sC_{sx} - \sum_{y<x} h_y(0) C_y$.
By construction $C_x$ satisfies the required conditions.
\end{proof}

\begin{lemma}\label{klunique}Suppose $C\in \sum_x v\ZZ[v]H_x$ is self dual. Then $C=0$.
\end{lemma}

\begin{proof}Write $C=\sum_x h_x H_x$ for $h_x\in v\ZZ[v,v^{-1}]$. Using \eqref{inverse} one sees that 
\[ \overline{H}_y \in H_y + \sum_{z<y} \ZZ[v, v^{-1}]H_z \]
for all $y\in W$. Pick $y$ maximal such that $h_y\neq 0$. Then $\overline{C}=C$ implies that $\overline{h}_y = h_y$. Thus, $h_y=0$, contradicting our assumption. So $C=0$. 
\end{proof}

\begin{thm}[\cite{KL79}]\label{klbasis}For each $x\in W$ there exists a unique self-dual element $C_x\in \cH$ such that $C_x \in H_x + \sum_{y} v\ZZ[v]H_y$.
\end{thm}

\begin{proof}The existence is given by Lemma \ref{klbasisstrong}, the uniqueness is provided by Lemma \ref{klunique}.
\end{proof}

\subsection{}For each $x,y\in W$ define polynomials $h_{y,x}\in \ZZ[v]$ by
\[ C_x = \sum_{y}h_{y,x} H_y.\]
The following result is well known. It will be key in the sequel.
\begin{lemma}[\cite{KL79}]\label{refineklalgorithm}Let $s\in S$ and let $x\in W$. If $sx<x$ then
\[ C_x = C_sC_{sx} - \sum_{\substack{y<sx, \\ sy<y}}\mu(y,sx)C_y \]
where $\mu(y,sx)$ is the coefficient of $v$ in $h_{y,sx}$.
\end{lemma}

\begin{proof}
We have
\begin{align*}
C_sC_{sx} &= C_s(H_{sx} + \sum_{y<sx}h_{y,sx}H_y \\
&= C_s\big(H_{sx} + \sum_{\substack{y<sx, \\ sy< y}} h_{y,sx}H_y + \sum_{\substack{y<sx, \\ sy>y}} h_{y,sx}H_y \big) \\
&=H_x + vH_{sx} + \sum_{\substack{y<sx, \\ sy< y}} h_{y,sx}(H_{sy}+v^{-1}H_y) + \sum_{\substack{y<sx, \\ sy>y}} h_{y,sx}(H_{sy}+vH_y).
\end{align*}
We need to look slightly more carefully at the construction of $C_x$ in Lemma \ref{klbasisstrong}.
Note that for $y<sx$ each $h_{y,sx}\in v\ZZ[v]$. So, by the proof of Lemma \ref{klbasisstrong},
\[ C_x = C_sC_{sx} - \sum_{\substack{y<sx, \\ sy<y}}\mu(y,sx)C_y. \qedhere\]
\end{proof}

\begin{prop}\label{throughswall}Let $s\in S$ and let $x\in W$. For each $y\in W$ let $\mu(y,x)$ denote the coefficient of $v$ in $h_{y,x}$, then
\[ C_sC_x= \begin{cases}
C_{sx} + \sum_{\substack{y<x, \\ sy < y}}\mu(y,x)C_y & \mbox{if $sx>x$}; \\
(v+v^{-1})C_x & \mbox{if $sx<x$}.\end{cases}\]
\end{prop}

\begin{proof}Suppose $sx<x$, then by Lemma \ref{refineklalgorithm}
\[ C_x = C_sC_{sx} - \sum_{\substack{y<sx, \\ sy<y}}\mu(y,sx)C_y. \]
Further, $C_s^2=(v+v^{-1})C_s$ by \eqref{inverse} and so the result follows by the induction hypothesis. If $sx>x$, then once again by Lemma \ref{refineklalgorithm}
\[ C_{sx} = C_sC_x - \sum_{\substack{y<x, \\ sy<y}}\mu(y,x)C_y,\]
whence the result.
\end{proof}

The following is originally due to Kazhdan-Lusztig \cite{KL79}, I learnt the proof presented here from \cite{So97}.
\begin{prop}[\cite{KL79}]\label{longest}Suppose $W$ is finite. Let $w_0\in W$ be the longest element. Then
\[ C_{w_0}=\sum_{x\in W} v^{\ell(w_0)-\ell(x)}H_x. \]
\end{prop}

\begin{proof}Let $H=\sum_y h_yH_y$, $h_y\in\ZZ[v,v^{-1}]$, be an element of $\cH$ such that $b(C_s)H = 0$ for all $s\in S$. Then, by Lemma \ref{translationthroughthewall}
\[
\sum_{sy>y} h_y(H_{sy}-v^{-1}H_y)+ \sum_{sy<y}h_y(H_{sy}-vH_y) =0 \]
for all $s\in S$. Consequently,
$h_s = v^{-1}h_e$
for all $s\in S$. Proceeding by induction we see that 
$H= h_e\sum_y v^{\ell(y)}H_y$.
By Prop.\ \ref{throughswall} $C_sC_{w_0} = (v+v^{-1})C_{w_0}$ for all $s\in S$. The result follows.
\end{proof}

\section{The Ext form}\label{s:ext}
\subsection{}Define a symmetric $\ZZ[v,v^{-1}]$-bilinear form $\langle \cdot, \cdot\rangle\colon \cH\otimes \cH \to \ZZ[v,v^{-1}]$ by
\begin{equation}
\langle H_x, H_y \rangle = \delta_{x,y}.
\end{equation}
\begin{remark}This form corresponds to the form considered in the proof of \cite[Thm.\ 3.11.4]{BGS}. Namely, it corresponds to
\[ \sum_{i,n} (-1)^i \dim\,\Ext^i(M, \DD N \langle n \rangle)v^{-n}, \]
where $M,N$ are perverse sheaves on the flag variety, $\DD$ is Verdier duality, $\langle n \rangle$ is shift + (half) Tate twist in the derived category and $v$ is a formal variable.
\end{remark}
\begin{prop}\label{adjointness}If $s\in S$, then
\[ \langle C_sH, H'\rangle = \langle H, C_sH' \rangle \]
for all $H,H'\in\cH$.
\end{prop}

\begin{proof}It suffices to prove the assertion for $H=H_x$ and $H'=H_y$ with $x,y\in W$. By Lemma \ref{translationthroughthewall} we have
\[ \langle C_sH_x, H_y \rangle = \begin{cases}
\delta_{sx, y} + v\delta_{x,y} & \mbox{if $sx>x$}; \\
\delta_{sx,y} + v^{-1}\delta_{x,y} & \mbox{if $sx<x$};
\end{cases}\]
and
\[ \langle H_x, C_sH_y \rangle = \begin{cases}
\delta_{x, sy} + v\delta_{x, y} & \mbox{if $sy>y$}; \\
\delta_{x,sy} + v^{-1}\delta_{x, y} & \mbox{if $sy<y$}.
\end{cases}\]
So there are four cases to consider: if $sx>x$ and $sy>y$, then the assertion is evident; if $sx>x$ and $sy<y$, then certainly $x\neq y$ and so the assertion holds; if $sx<x$ and $sy>y$, then once again $x\neq y$ and the assertion holds; finally, if $sx<x$ and $sy<y$, then the assertion is evident.
\end{proof}

\subsection{}Define a ring involution $b\colon \cH\to \cH$ by
\begin{equation}\label{defb} b(v)=-v^{-1},\qquad b(H_x)=H_x.\end{equation}
Then $b$ commutes with $d$. Furthermore:
\begin{prop}\label{bcadjointness}If $s\in S$, then
\[ \langle b(C_s)H, H'\rangle = \langle H, b(C_s)H'\rangle \]
for all $H,H'\in\cH$.
\end{prop}

\begin{proof}Since $b(C_s) = C_s-(v+v^{-1})$, the result follows from Prop.\ \ref{adjointness}.
\end{proof}

\begin{lemma}\label{projectivebasis}Assume $W$ is finite. Then for each $x\in W$ there exists a $P_x\in H_x + \sum_{y>x}v\ZZ[v]H_y$ such that
$\langle P_x, C_z \rangle = \delta_{x,z}$
for all $z\in W$.
\end{lemma}

\begin{proof}Proceed by induction. Let $w_0$ be the longest element in $W$. We start our induction with $P_{w_0} = H_{w_0}$. Now let $x\in W$ be given and suppose we know the existence of $P_y$ for all $y>x$. If $x\neq w_0$ find $s\in S$ such that $sx>x$. Then by the induction hypothesis
\[ b(C_s)P_{sx} \in H_x + \sum_{y>x}h_yH_y\]
for some $h_y\in \ZZ[v]$.
For each $y > x$ let $p_y = \langle b(C_s)P_{sx}, C_y\rangle$. Set
\[ P_x = b(C_s)P_{sx} - \sum_{y>x} p_yP_y.\]
Let's show that $P_x \in H_x + \sum_{y>x}v\ZZ[v]H_y$. To do this we must demonstrate that (a) $p_y\in \ZZ[v]$ and (b) $p_y(0)=h_y(0)$, for each $y>x$. Both of these follow from the fact that each $h_y\in\ZZ[v]$ and each $C_y \in H_y + \sum_{z}v\ZZ[v]H_z$.
By construction $\langle P_x, C_z \rangle = \delta_{x,z}$ for all $z\in W$.
\end{proof}

\begin{thm}\label{projectivesexist}Assume $W$ is finite. Then for each $x\in W$ there exists a unique $P_x\in\cH$ such that
$\langle P_x, C_y \rangle =\delta_{x,y}$
for all $y\in W$.
\end{thm}

\begin{proof}Existence was established in the previous Lemma. Uniqueness follows from the evident fact that the form $\langle \cdot, \cdot \rangle$ is non-degenerate.
\end{proof}

\section{Dualities between the bases}\label{s:main}
\begin{lemma}\label{keylemma}Assume $W$ is finite. Let $w_0$ be the longest element in $W$. Then for all $x\in W$,
$P_x = C^xH_{w_0}$
for some self-dual $C^x\in\cH$.
\end{lemma}

\begin{proof}Proceed by induction. We start our induction with $P_{w_0}=H_{w_0}$. Now let $x\in W$ be given and suppose we know that $P_y = C^yH_{w_0}$, with $\overline{C^y}=C^y$, for all $y>x$. If $x\neq w_0$ find $s\in S$ such that $sx>x$. Then, by the proof of Lemma \ref{projectivebasis} and the induction hypothesis,
\[
P_x = b(C_s)P_{sx} - \sum_{y>x}p_yP_y 
= b(C_s)C^{sx}H_{w_0} - \sum_{y>x}p_yC^yH_{w_0}, \]
where
$p_y = \langle b(C_s)P_{sx}, C_y \rangle$.
To complete the proof it suffices to show that $p_y\in \ZZ[v+v^{-1}]$ for each $y$. By Prop.\ \ref{bcadjointness}
\[ \langle b(C_s)P_{sx}, C_y\rangle = \langle P_{sx}, b(C_s)C_y \rangle.\]
If $sy<y$, then $b(C_s)C_y=(C_s-(v+v^{-1}))C_y=0$ by Prop.\ \ref{throughswall}. If $sy>y$, then once again using Prop.\ \ref{throughswall},
\[
b(C_s)C_y = (C_s-(v+v^{-1}))C_y 
= C_{sy}- (v+v^{-1})C_y +\sum_z m_z C_z
\]
for appropriate integers $m_z\in\ZZ$. The result follows.
\end{proof}

\begin{thm}\label{tiltingduality}Assume $W$ is finite. Let $w_0$ be the longest element in $W$. Then
\[ b(C_x)H_{w_0} = P_{xw_0} \]
for all $x\in W$.
\end{thm}

\begin{proof}
By Lemma \ref{keylemma} $P_{xw_0}H_{w_0}^{-1}$ is self-dual. On the other hand, by Lemma \ref{projectivebasis},
\[ P_{xw_0}H_{w_0}^{-1} \in H_{xw_0}H_{w_0}^{-1} +\sum_{y>xw_0} v\ZZ[v]H_yH_{w_0}^{-1}. \]
Using the fact that $H_zH_{w_0}^{-1} = H_zH_z^{-1} H_{w_0z^{-1}}^{-1} = H_{w_0z^{-1}}^{-1}$ and applying $d$ we obtain
\[ P_{xw_0}H_{w_0}^{-1} \in H_{x} + \sum_{y>xw_0}v^{-1}\ZZ[v^{-1}]H_{yw_0}.\]
By the uniqueness of $C_{xw_0}$ (Thm.\ \ref{klbasis}) we must have that 
\[P_{xw_0}H_{w_0}^{-1} = b(C_x).\qedhere\]
\end{proof}

\end{document}